\documentclass[12pt]{article}
\usepackage{amsmath}
\usepackage{color}
\usepackage{fancyhdr}
\usepackage{verbatim}
\usepackage{amsfonts,longtable,mathrsfs}
\usepackage{epsfig}
\usepackage{amsfonts}
\usepackage{color}
\usepackage[numbers,sort&compress]{natbib}

\def\R{\hbox{{\rm I}\kern-0.2em{\rm R}\kern0.2em}}

\def\bn{\begin{equation}}
\def\en{\end{equation}}
\def\bny{\begin{eqnarray}}
\def\eny{\end{eqnarray}}
\def\be{\begin{eqnarray*}}
\def\ee{\end{eqnarray*}}
\def\bc{\begin{center}}
\def\ec{\end{center}}

\def\({\left(}
\def\){\right  )}
\def\[{\left[}
\def\]{\right]}
\def\bc{\begin{center}}
\def\ec{\end{center}}

\newtheorem{dfn}{Definition}[section]
\newtheorem{thm}{Theorem}[section]
\newtheorem{rem}{Remark}[section]
\newtheorem{pro}{Proposition}[section]

\newtheorem{cor}{Corollary}[section]
\newtheorem{lem}{Lemma}[section]
\newtheorem{exm}{Example}[section]

\def\bn{\begin{equation}}
\def\en{\end{equation}}
\def\bny{\begin{eqnarray}}
\def\eny{\end{eqnarray}}
\def\be{\begin{eqnarray*}}
\def\ee{\end{eqnarray*}}
\def\bdn{\begin{dfn}}
\def\edn{\end{dfn}}
\def\btm{\begin{thm}}
\def\etm{\end{thm}}
\def\bpf{\begin{proof}}
\def\epf{\end{proof}}
\def\bpn{\begin{pro}}
\def\epn{\end{pro}}
\def\brk{\begin{rem}}
\def\erk{\end{rem}}
\def\bcy{\begin{cor}}
\def\ecy{\end{cor}}
\def\blm{\begin{lem}}\def\elm{\end{lem}}
\def\bex{\begin{exm}}
\def\eex{\end{exm}}

 \def\R{{\hat R}}
\begin{document}

\bc {\bf On certain sixth order difference equations
  }\ec
\medskip
\bc
 D. Nyirenda$^{1,}$
and  M. Folly-Gbetoula$^{1,}$
 \vspace{1cm}
\\$^{1}$ School of Mathematics, University of the Witwatersrand, Johannesburg, South Africa.\\

\ec
\begin{abstract}
\noindent We use the Lie group analysis method to investigate the invariance properties and the solutions of
\begin{align*}
x_{n+1} =\frac{x_{n-5}x_{n-3}}{x_{n-1}(a_n +b_nx_{n-5}x_{n-3})}.
\end{align*}
We show that this equation has a two-dimensional Lie algebra and that its solutions can be presented in a unified manner. Besides presenting solutions of the recursive sequence above where $(a_n)$ and $(b_n)$ are sequences of real numbers, some specific cases are emphasized.
\end{abstract}
\section{Introduction} \setcounter{equation}{0}
Difference equations are important in mathematical modeling, especially where discrete time evolving variables are concerned.
They also occur when studying discrimination methods for differential equations. Countless results in the subject of difference equations have been recorded \cite{bc,ibm0, ibm1,ibm2,KE,IY, El1,El2}. For rational difference equations of order greater than two, the study  can be quite challenging at the same time rewarding. Rewarding in the sense that such a study lays ground for the theory of global properties of difference equations (not necessarily rational) of higher-order.
\par \noindent In \cite{hydon0}, Hydon developed an effective symmetry based algorithm to deal with the obtention of solutions of difference equations of any order. However, the calculation one deals with in this application to difference equations of order greater than one can become cumbersome. The method consists of finding a group of transformations that maps solutions onto themselves. Symmetry method is a valuable tool and it has been used to solve several difference equations \cite{FK, FK2, MNF, DM}.
 \par \noindent
 In this paper, our objective is to obtain the symmetry operators of
\begin{align}\label{xn}
x_{n+1} =\frac{x_{n-5}x_{n-3}}{x_{n-1}(a_n +b_nx_{n-5}x_{n-3})},
\end{align}
where ${(a_n)}_{\mathbb{N}_0}$ and ${(b_n)}_{\mathbb{N}_0}$ are real sequences and to find its solutions by way of symmetries. Without loss of generality, we equivalently study the forward difference equation
\begin{align}\label{1.0}
u_{n+6} =\frac{u_nu_{n+2}}{u_{n+4}(A_n +B_nu_nu_{n+2})}.
\end{align}
We refer the interested reader to \cite{hydon0, PO} for a deeper knowledge of Lie analysis.
\section{Definitions and Notation}
In this section, we briefly present some definitions and notation (largely from Hydon in \cite{hydon0}) indispensable for the understanding of Lie symmetry analysis of difference equations.
\begin{dfn}
 Let $G$ be a local group of transformations acting on a manifold $M$. A subset $\mathcal{S}\subset M$ is called G-invariant, and $G$ is called symmetry group of $\mathcal{S}$, if whenever $x\in \mathcal{S} $, and $g\in G$ is such that $g\cdot x$ is defined, then $g\cdot x \in \mathcal{S}$.
\end{dfn}
\begin{dfn}
 Let $G$ be a connected group of transformations acting on a manifold $M$. A smooth real-valued function $\mathcal{V}: M\rightarrow \mathbb{R}$ is an invariant function for $G$ if and only if $$X(\mathcal{V})=0\qquad \text { for all } \qquad  x\in M,$$
and every infinitesimal generator $X$ of $G$.
\end{dfn}
\begin{dfn}
A parameterized set of point transformations,
\begin{equation}
\Gamma_{\varepsilon} :x\mapsto \hat{x}(x;\varepsilon),
\label{eq: b}
\end{equation}
where $x=x_i, $ $i=1,\dots,p$ are continuous variables, is a one-parameter local Lie group of transformations if the following conditions are satisfied:
\begin{enumerate}
\item $\Gamma_0$ is the identity map if $\hat{x}=x$ when $\varepsilon=0$
\item $\Gamma_a\Gamma_b=\Gamma_{a+b}$ for every $a$ and $b$ sufficiently close to 0
\item Each $\hat{x_i}$ can be represented as a Taylor series (in a neighborhood of $\varepsilon=0$ that is determined by $x$), and therefore
\end{enumerate}
\begin{equation}
\hat{x_i}(x:\varepsilon)=x_i+\varepsilon \xi _i(x)+O(\varepsilon ^2), i=1,...,p.
\label{eq: c}
\end{equation}
\end{dfn}
Assuming that the sixth-order difference equation has the form
\begin{align}\label{general}
u_{n+r}=&\Psi (n, u_n,  \dots, u_{n+r-1}), \quad n\in D
\end{align}
for some smooth function $\Omega$ and a regular domain $D\subset \mathbb{Z}$.
To deduce the symmetry group of \eqref{general}, we search for a one parameter Lie group of point transformations
\begin{equation}\label{Gtransfo}
\Gamma_{\varepsilon}: (n,u_n) \mapsto(n, u_n+\varepsilon Q (n,u_n)),
\end{equation}
in which $\varepsilon$ is the parameter and $Q$ a continuous function, referred to as characteristic. Let
\begin{align}\label{Ngener}
X= & Q(n,u_n)\frac{\partial}{ \partial u_n}+SQ(n,u_{n})\frac{\partial\quad }{ \partial u_{n+1}} +\dots+ S^{r-1}Q(n,u_{n})\frac{\partial\quad }{ \partial u_{n+r-1}}
\end{align}
be the corresponding $\lq$prolonged\rq \, infinitesimal generator and  $S: n\mapsto n+1$ the shift operator. The linearized symmetry condition is given by
\begin{align}\label{LSC}
 S^r Q- X \Psi =0.
\end{align}
Upon knowledge of the characteristic $Q$, it is important to introduce the canonical coordinate 
 \begin{align}\label{cano}
S_n= \int{\frac{du_n}{Q(n,u_n)}},
 \end{align}
a useful tool which allows one to obtain the invariant function of $\Gamma _{\varepsilon}$.
\section{Main results}
As earlier emphasized, our equation under study is
\begin{align}\label{un}
u_{n+6}=\Psi =\frac{u_nu_{n+2}}{u_{n+4}(A_n +B_nu_nu_{n+2})}.
\end{align}
Appliying the criterion of invariance \eqref{LSC} to \eqref{un}, we get
\begin{align}\label{a1}
 &Q(n+6,\Psi)+
  {\frac{ {u_n}u_{n+2}Q\left({n+4, u_{n+4}}\right)}{{u_{n+4}^2\left(A_n +B_nu_nu_{n+2}\right)}} }-\frac{A_nu_nQ\left(n+2,{u_{n+2}}\right) }{{u_{n+4}\left((A_n +B_nu_nu_{n+2}\right)}^{2}} \nonumber\\&-{\frac{A_n u_{n+2}Q\left({n, u_{n}}\right)}{{u_{n+4}\left(A_n +B_nu_nu_{n+2}\right)}^{2}} } =0.
\end{align}
In order to eliminate $\Psi$ in \eqref{a1}, we invoke implicit differentiation with respect to $u_n$ (regarding $u_{n+4}$ as a function of $u_n$, $u_{n+2}$ and $\Psi$) via the operator
 $$L=\frac{\partial\quad }{\partial {u_n}}+\left(\frac{\partial{u_{n+4}}}{\partial{u_{n}\quad}}\right)\frac{\partial\qquad }{\partial u_{n+4}}=\frac{\partial\quad }{\partial {u_n}}-\left(\frac{\Psi _{,u_n}}{\Psi _{,u_{n+4}}}\right)\frac{\partial\qquad }{\partial u_{n+4}}.$$
 Note that $f_{,x}$ denotes the derivative with respect to $x$.
With some simplification, one gets
\begin{align}\label{a3}
& \left(A_n +B_nu_nu_{n+2}\right) Q'\left(n+4, u_{n+4}\right) - \frac{\left(A_n +B_nu_nu_{n+2}\right)Q\left(n+4, u_{n+4}\right)}{ u_{n+4}}\nonumber\\
 &+ B_n u_n Q(n+2,u_{n+2})-\left(A_n +B_nu_nu_{n+2}\right)Q'(n,u_n)\nonumber\\
 &+\left(2B_n u_{n+2}+\frac{A_n}{ u_{n}}\right)Q\left(n, u_{n}\right) =0.
\end{align}
Note that the symbol $'$ stands for the derivative with respect to the continuous variable. After twice differentiating (\ref{a3}) with respect to $u_{n}$, keeping $u_{n+2}$ and $u_{n+4}$ fixed, we are led to the equation
\begin{align}\label{a4'}
&-B_n u_n u_{n+2}Q'''(n,u_n)-A_nQ'''(n,u_n)+\frac{A_n}{u_n}
Q''\left(n,{u_{n}}\right)-\frac{2A_n}{{u_n}^2}Q'\left(n,u_{n}\right)\nonumber\\
 &+ \frac{2A_n}{{u_n}^3}Q\left(n,u_{n}\right)=0.
\end{align}
\noindent Note that the characteristic in \eqref{a4'} is not a function of $u_{n+2}$ and so we split \eqref{a4'} up with respect to powers of $u_{n+2}$ to get the system
\begin{subequations}\label{a4}
\begin{align}
 1&:Q'''(n,u_n)-\frac{1}{u_n}
Q''\left(n,{u_{n}}\right)+\frac{2}{{u_n}^2}Q'\left(n,u_{n}\right)
 - \frac{2}{{u_n}^3}Q\left(n,u_{n}\right)=0\\&\nonumber\\
 u_{n+2}&:Q'''(n,u_n)=0.
 \end{align}
\end{subequations}
\noindent We find that the solution to \eqref{a4} is
\begin{align}\label{a6}
Q\left(n,{u_{n}}\right) = \alpha_ n  {u_n}^2 +\beta _n  {u_n}
\end{align}
for some arbitrary functions $\alpha _n$ and $\beta _n$ that depend on $n$. Substituting  \eqref{a6} and its shifts  in \eqref{a1}, and then replacing the expression of $u_{n+6}$ given in \eqref{un} in the resulting equation yields
\begin{align}\label{a7}
& B_n{u_{n}}^2{u_{n+2}}^2{u_{n+4}}^2\alpha_{n+4}+
B_n{u_{n}}^2{u_{n+2}}^2{u_{n+4}}(\beta_{n+4}+\beta_{n+3})\nonumber \\
&
-A_n {u_{n}}^2{u_{n+2}}{u_{n+4}}\alpha_{n}
-A_n {u_{n}}{u_{n+2}}^2{u_{n+4}}\alpha_{n+2}
+A_n {u_{n}}{u_{n+2}}{u_{n+4}}^2\alpha_{n+4}\nonumber\\
&
+{u_{n}}^2{u_{n+2}}^2\alpha_{n+2}
-A_n\left(  \beta_{n} +\beta_{n+2}-\beta_{n+4}-\beta_{n+3} \right)=0.
\end{align}
Equating all coefficients of all powers of shifts of $u_n$ to zero and simplifying the resulting system, we get its reduced form
\begin{align}
&\alpha _n=0,\label{rel1}  \\
&\beta_{n} + \beta _{n+2} =0.\label{rel2}
\end{align}
The two independent solutions of the linear second-order difference equation \eqref{rel2} are given by
\begin{align}
\beta_n=\beta ^n  \text{ and  } \beta_n=\bar{\beta}^n,
\end{align}
where $\beta = \exp \{i\pi/2\}$ and $\bar{\beta}=-\exp \{i\pi/2\}$ is its complex conjugate.
The characteristic functions are given by
\begin{align}\label{a9}
Q_1(n,u_n)=\beta ^n u_n \quad \text{ and  }\quad  Q_2(n,u_n)=\bar{\beta}^n u_n,
\end{align}
and so the symmetry operators are
\begin{subequations}\label{gener}
\begin{align}
X_1= &\beta ^{n} u_{n} \frac{\partial}{\partial { u_{n}}}+\beta ^{n+1} u_{n+1} \frac{\partial}{\partial { u_{n+1}}}
+\beta ^{n+2} u_{n+2} \frac{\partial}{\partial { u_{n+2}}}+\beta ^{n+3} u_{n+3} \frac{\partial}{\partial { u_{n+3}}}\nonumber\\
&+\beta ^{n+4} u_{n+4} \frac{\partial}{\partial { u_{n+4}}}+\beta ^{n+5} u_{n+5} \frac{\partial}{\partial { u_{n+5}}},\\& \nonumber \\
X_2= & \bar{\beta} ^{n} u_{n} \frac{\partial}{\partial { u_{n}}}+\bar{\beta} ^{n+1} u_{n+1} \frac{\partial}{\partial { u_{n+1}}}
+\bar{\beta} ^{n+2} u_{n+2} \frac{\partial}{\partial { u_{n+2}}}+\bar{\beta} ^{n+3} u_{n+3} \frac{\partial}{\partial { u_{n+3}}}\nonumber\\
&+\bar{\beta} ^{n+4} u_{n+4} \frac{\partial}{\partial { u_{n+4}}}+\bar{\beta} ^{n+5} u_{n+5} \frac{\partial}{\partial { u_{n+5}}}.
\end{align}
\end{subequations}
Using the canonical coordinate
\begin{align}\label{cano}
S_n =\int\frac{du_n}{Q_1(n,u_n)}=\int\frac{du_n}{\beta ^nu_n}=\frac{1}{\beta ^n}\ln|u_n|
\end{align}
and \eqref{rel2}, we derive the invariant function $\tilde{\mathcal{V}}_n$ as follows:
\begin{align}\label{tilde}
\tilde{\mathcal{V}}_n = S_n\beta^ n + S_{n+2}\beta^{n+2}.
\end{align}
Actually,
\begin{align}
X_1 (\tilde{\mathcal{V}}_n) = \beta^ n+ \beta^{n+2}=0
\end{align}
and
\begin{align}
X_2 (\tilde{\mathcal{V}}_n) = \bar{\beta}^ n + \bar{\beta}^{n+2}=0.
\end{align}
For the sake of convenience, we use
\begin{align}\label{vn}
|{\mathcal{V}}_n| =\exp\{  -\tilde{\mathcal{V}}_n \},
\end{align}
instead. In other words, $\mathcal{V}_n =\pm 1/(u_nu_{n+2})$.
Using \eqref{un} and \eqref{vn}, one can prove that
\begin{equation}\label{vn1}
\mathcal{V}_{n+4}={A_n\mathcal{V}_n}\pm B_n.
\end{equation}
From here, to obtain the solution of \eqref{un}, we first employ \eqref{cano} to get
\begin{align}
\vert u_{n}\vert  =& \exp\left(\beta _n S_n\right).
\end{align}
Secondly, we employ \eqref{tilde} to obtain
\begin{align}
\vert u_{n}\vert
=&\exp\left(\beta^{n}c_1 + \bar{\beta}^{n}c_2 - \frac{1}{2}\sum_{k_1 = 0}^{n - 1}{\beta}^{n}\bar{\beta}^{k_1}\tilde{\mathcal{V}}_{k_1} - \frac{1}{2}\sum_{k_2 = 0}^{n - 1} \bar{\beta}^{n}{\beta}^{k_2}\tilde{\mathcal{V}}_{k_2}\right).
\end{align}
Lastly, we use \eqref{vn} to get
\begin{subequations}\label{solsingle}
\begin{align}\label{single}
\vert u_{n}\vert
           =& \exp\left(\beta^{n}c_1 + \bar{\beta}^{n}c_2 + \frac{1}{2}\sum\limits_{k_1 = 0}^{n - 1}{\beta}^{n}\bar{\beta}^{k_1}\ln|{\mathcal{V}}_{k_1}| + \frac{1}{2}\sum\limits_{k_2 = 0}^{n - 1} \bar{\beta}^{n}{\beta}^{k_2}\ln|{\mathcal{V}}_{k_2}|\right),\nonumber\\
           =&\exp\left(\beta^{n}c_1 + \bar{\beta}^{n}c_2 + \sum\limits_{k = 0}^{n - 1}\text{Re}[\gamma(n,k)]\ln|{\mathcal{V}}_{k}| \right)
\end{align}
in which $V_k$ satisfies \eqref{vn1} with $\gamma(n,k)=\beta ^n \bar{\beta}^{k}$. Note that the constants $c_{1}$ and $c_{2}$ satisfy
 \begin{align}
 c_1+c_2&=\ln |u_0| \quad \text{ and } \quad
 \beta( c_{1} - c_{2})=\ln |u_1|.
\end{align}
\end{subequations}
It is worthwhile to mention that equations in (\ref{solsingle}) give the solutions of (\ref{1.0}) in a unified manner.\par \noindent
On a further note, $\gamma(n, k)=\beta ^n \bar{\beta}^{k}$ satisfies
\begin{align}\label{m11'}
&\gamma (0,1)=\bar{\beta}, \gamma (1,0)={\beta}, \gamma (n,n)=1, \gamma (n+2,k)=-\gamma (n,k),\nonumber\\
&\gamma (n,k+2)=-\gamma (n,k),\gamma (4n,k)=\gamma(0,k), \gamma(n,4k)=\gamma({n,0}).
\end{align}
From $u_n$ given in \eqref{single} and properties \eqref{m11'}, observe that
\begin{equation}\label{diff1}
\vert u_{4n + j}\vert  = \exp\left(H_{j} +  \sum\limits_{k = 0}^{4n + j - 1}\text{Re}(\gamma(j,k))\ln\vert \mathcal{V}_{k} \vert \right),
\end{equation}
$j=0,1,2,3$, in which $$H_{j} = \beta^{j}c_1 + \bar{\beta}^{j}c_2.$$
\noindent For $j = 0$, we have,
\begin{align}\label{4n}
 \vert u_{4n} \vert = &\exp(H_{0} + \ln \vert \mathcal{V}_{0}\vert - \ln \vert \mathcal{V}_{2}\vert
   + \ldots +  \ln \vert \mathcal{V}_{4n - 4}\vert - \ln \vert \mathcal{V}_{4n - 2}\vert ) \nonumber\\= & \exp(H_{0})\prod_{s = 0}^{n - 1}\left| \frac{\mathcal{V}_{4s}}{\mathcal{V}_{4s + 2}}\right|. \end{align}
In order to deduce $\exp(H_{0})$, we set $n = 0$ in \eqref{4n} and note that $\vert u_0 \vert = \exp(H_{0})$. \par \noindent
It can be shown that there is no need for the absolute values via the utilization of the fact that
\begin{equation}\label{eqq}
\mathcal{V}_{i} = \frac{1}{u_{i}u_{i+ 2}}.
\end{equation}
Thus
 $$ u_{4n} = u_0\prod_{s = 0}^{n - 1}\frac{\mathcal{V}_{4s}}{\mathcal{V}_{4s + 2}}.$$
 Similarly, for $j = 0,1,2,3$, we obtain
\begin{equation}\label{eqv1}
u_{4n + j} = u_j\prod\limits_{s = 0}^{n - 1}\frac{\mathcal{V}_{4s + j}}{\mathcal{V}_{4s + j + 2}}.
\end{equation}
Nevertheless, from \eqref{vn1}, using the plus sign we are led to
\begin{align}\label{eqv2}
\mathcal{V}_{4n+j}\quad=&\mathcal{\mathcal{V}}_j \left(   \prod\limits_{k_1=0}^{n-1}A_{4k_1+j}\right) +\sum\limits_{l=0}^{n-1} \left(  B_{4l+j}\prod\limits_{k_2=l+1}^{n-1}A_{4k_2+j}\right),
\end{align}
for $\; j=0,1,2,3,$ where $\mathcal{V}_{j} = {1}/({u_ju_{j+2}})$. Thus, using \eqref{eqv1} and \eqref{eqv2} with $j = 0$, we get
\begin{align}
u_{4n} & = u_{0}\prod\limits_{s = 0}^{n - 1}\frac{\mathcal{V}_{4s}}{\mathcal{V}_{4s + 2}}\nonumber\\
       & = u_0\prod\limits_{s = 0}^{n - 1}\frac{\mathcal{V}_{0}\left(\prod\limits_{k_1 = 0}^{s - 1}A_{4k_1}\right) + \sum\limits_{l = 0}^{s - 1}\left(B_{4l}\prod\limits_{k_2 = l + 1}^{s - 1}A_{4k_2}\right) }{ \mathcal{V}_{2}\left(\prod\limits_{k_1 = 0}^{s - 1}A_{4k_1 + 2}\right) + \sum\limits_{l = 0}^{s - 1}\left(B_{4l + 2}\prod\limits_{k_2 = l + 1}^{s - 1}A_{4k_2 + 2}\right) }\nonumber\\
       & = u_0\prod\limits_{s = 0}^{n - 1}\frac{u_4}{u_0}\frac{\left(\prod\limits_{k_1 = 0}^{s - 1}A_{4k_1}\right) + u_0u_2\sum\limits_{l = 0}^{s - 1}\left(B_{4l}\prod\limits_{k_2 = l + 1}^{s - 1}A_{4k_2}\right) }{\left(\prod\limits_{k_1 = 0}^{s - 1}A_{4k_1 + 2}\right) + u_2u_4\sum\limits_{l = 0}^{s - 1}\left(B_{4l + 2}\prod\limits_{k_2 = l + 1}^{s - 1}A_{4k_2 + 2}\right) }\nonumber\\
       & = u_0^{1 - n}u_4^{n}\prod\limits_{s = 0}^{n - 1}\frac{\left(\prod\limits_{k_1 = 0}^{s - 1}A_{4k_1}\right) + u_0u_2\sum\limits_{l = 0}^{s - 1}\left(B_{4l}\prod\limits_{k_2 = l + 1}^{s - 1}A_{4k_2}\right) }{\left(\prod\limits_{k_1 = 0}^{s - 1}A_{4k_1 + 2}\right) + u_2u_4\sum\limits_{l = 0}^{s - 1}\left(B_{4l + 2}\prod\limits_{k_2 = l + 1}^{s - 1}A_{4k_2 + 2}\right) }.
\end{align}
\noindent For $j = 1$, we have
\begin{align}
u_{4n + 1} & = u_1\prod\limits_{s = 0}^{n - 1}\frac{\mathcal{V}_{4s + 1}}{\mathcal{V}_{4s + 3}}\nonumber\\
           & = u_1^{1 - n}u_5^{n}\prod\limits_{s = 0}^{n - 1}\frac{\left(\prod\limits_{k_1 = 0}^{s - 1}A_{4k_1 + 1}\right) + u_1u_3\sum\limits_{l = 0}^{s - 1}\left(B_{4l + 1}\prod\limits_{k_2 = l + 1}^{s - 1}A_{4k_2 + 1}\right)}{\left(\prod\limits_{k_1 = 0}^{s - 1}A_{4k_1 + 3}\right) + u_{3}u_{5}\sum\limits_{l = 0}^{s - 1}\left(B_{4l + 3}\prod\limits_{k_2 = l + 1}^{s - 1}A_{4k_2 + 3}\right)}.
\end{align}
For $j = 2$, we have
\begin{align}
u_{4n + 2} & = u_2\prod\limits_{s = 0}^{n - 1}\frac{\mathcal{V}_{4s + 2}}{\mathcal{V}_{4s + 4}}\nonumber\\
           & = u_0^{n}u_4^{-n}u_2\prod\limits_{s = 0}^{n - 1}\frac{ \left(\prod\limits_{k_1 = 0}^{s - 1}A_{4k_1 + 2}\right) + u_2u_4\sum\limits_{l = 0}^{s - 1}\left(B_{4l + 2}\prod\limits_{k_2 = l + 1}^{s - 1}A_{4k_2 + 2}\right)}{\left(\prod\limits_{k_1 = 0}^{s}A_{4k_1}\right) + u_0u_2\sum\limits_{l = 0}^{s}\left(B_{4l}\prod\limits_{k_2 = l + 1}^{s}A_{4k_2}\right)}.
\end{align}
For $j = 3$, we get
\begin{align}
u_{4n + 3} & = u_3\prod\limits_{s = 0}^{n - 1}\frac{\mathcal{V}_{4s + 3}}{\mathcal{V}_{4s + 5}} \nonumber\\
           & = u_1^{n}u_5^{-n}u_3\prod\limits_{s = 0}^{n - 1}\frac{\left(\prod\limits_{k_1 = 0}^{s - 1}A_{4k_1 + 3}\right) + u_3u_5\sum\limits_{l = 0}^{s - 1}\left(B_{4l + 3}\prod\limits_{k_2 = l + 1}^{s - 1}A_{4k_2 + 3}\right)}{\left(\prod\limits_{k_1 = 0}^{s}A_{4k_1 + 1}\right) + u_1u_3\sum\limits_{l = 0}^{s}\left(B_{4l + 1}\prod\limits_{k_2 = l + 1}^{s}A_{4k_2 + 1}\right)}.
\end{align}
Hence, 
 our solution in terms of $x_{n}$ is given by
\begin{equation}\label{eqn0}
x_{4n - 5} = \frac{x_{-1}^{n}}{x_{-5}^{n-1}}\prod\limits_{s = 0}^{n - 1}\frac{\left(\prod\limits_{k_1 = 0}^{s - 1}a_{4k_1}\right) + x_{-5}x_{-3}\sum\limits_{l = 0}^{s - 1}\left(b_{4l}\prod\limits_{k_2 = l + 1}^{s - 1}a_{4k_2}\right) }{\left(\prod\limits_{k_1 = 0}^{s - 1}a_{4k_1 + 2}\right) + x_{-3}x_{-1}\sum\limits_{l = 0}^{s - 1}\left(b_{4l + 2}\prod\limits_{k_2 = l + 1}^{s - 1}a_{4k_2 + 2}\right) },
\end{equation}
\begin{equation}\label{eqn1}
x_{4n - 4} =  \frac{x_{0}^{n} }{x_{-4}^{ n-1}}\prod\limits_{s = 0}^{n - 1}\frac{\left(\prod\limits_{k_1 = 0}^{s - 1}a_{4k_1 + 1}\right) + x_{-4}x_{-2}\sum\limits_{l = 0}^{s - 1}\left(b_{4l + 1}\prod\limits_{k_2 = l + 1}^{s - 1}a_{4k_2 + 1}\right)}{\left(\prod\limits_{k_1 = 0}^{s - 1}a_{4k_1 + 3}\right) + x_{-2}x_{0}\sum\limits_{l = 0}^{s - 1}\left(b_{4l + 3}\prod\limits_{k_2 = l + 1}^{s - 1}a_{4k_2 + 3}\right)},
\end{equation}
\begin{equation}\label{eqn2}
x_{4n - 3} = \frac{x_{-5}^{n}x_{-3}}{x_{-1}^{n}}\prod\limits_{s = 0}^{n - 1}\frac{ \left(\prod\limits_{k_1 = 0}^{s - 1}a_{4k_1 + 2}\right) + x_{-3}x_{-1}\sum\limits_{l = 0}^{s - 1}\left(b_{4l + 2}\prod\limits_{k_2 = l + 1}^{s - 1}a_{4k_2 + 2}\right)}{\left(\prod\limits_{k_1 = 0}^{s}a_{4k_1}\right) + x_{-5}x_{-3}\sum\limits_{l = 0}^{s}\left(b_{4l}\prod\limits_{k_2 = l + 1}^{s}a_{4k_2}\right)}
\end{equation}
and
\begin{equation}\label{eqn3}
x_{4n - 2} =\frac{ x_{-4}^{n}x_{-2}}{x_{0}^{n}}\prod\limits_{s = 0}^{n - 1}\frac{\left(\prod\limits_{k_1 = 0}^{s - 1}a_{4k_1 + 3}\right) + x_{-2}x_{0}\sum\limits_{l = 0}^{s - 1}\left(b_{4l + 3}\prod\limits_{k_2 = l + 1}^{s - 1}a_{4k_2 + 3}\right)}{\left(\prod\limits_{k_1 = 0}^{s}a_{4k_1 + 1}\right) + x_{-4}x_{-2}\sum\limits_{l = 0}^{s}\left(b_{4l + 1}\prod\limits_{k_2 = l + 1}^{s}a_{4k_2 + 1}\right)}
\end{equation}
with the conditions for well definedness given by
\begin{equation}\label{welldef}
- x_{-5+j}x_{-3+j}\sum\limits_{l = 0}^{s-i}\left(b_{4l + j}\prod\limits_{k_2 = l + 1}^{s-i}a_{4k_2 + j}\right)\neq \left(\prod\limits_{k_1 = 0}^{s-i}a_{4k_1 + j}\right),
\end{equation}
for $i=0, \;j\in \{0,1\}$; and $i=1,\;j\in\{2,3\}$.
\vspace{.5cm}

\par \noindent
In the following sections, we let  $x_{-5}=c,\;x_{-4}=d,\;x_{-3}=e,\;x_{-2}=f,\;x_{-1}=g,\;x_0=h$ and we specifically look at some special cases.
\section{The case $a_{n}$ and $b_{n}$ are 1-periodic}
Let $a_n = a$ and $b_n = b$, where $a$ and $b$ are non-zero constants.
\subsection{Case: $a \neq 1$}
We have\\
\begin{equation*}\label{sole0}
x_{4n - 5} = \frac{g^{n}}{c^{ n-1}}\prod_{s = 0}^{n - 1}\frac{ a^{s} + bce \frac{1 - a^{s}}{1 - a} }{ a^{s} + beg \frac{1 - a^{s}}{1 - a}},\;
 x_{4n - 4} =\frac{ h^{n}}{d^{ n-1}}\prod_{s = 0}^{n - 1}\frac{a^{s} + bdf\frac{1 - a^{s}}{1 - a}}{a^{s} + bfh\frac{1 - a^{s}}{1 - a}},
\end{equation*}
\begin{equation*}\label{sole2}
x_{4n - 3} = \frac{c^{n}e}{g^{n}}\prod_{s = 0}^{n - 1}\frac{ a^{s} + beg\frac{1 - a^{s}}{1 - a}}{a^{s+1} + bce \frac{1 - a^{s + 1}}{1 - a}} ,\;
x_{4n - 2} = \frac{d^{n}f}{h^{n}}\prod_{s = 0}^{n - 1}\frac{a^{s} + bfh\frac{1 - a^{s}}{1 - a}}{a^{s + 1} + bdf \frac{1 - a^{s + 1}}{1 - a}}.
\end{equation*}
Here, condition \eqref{welldef} becomes
\begin{align}
-x_{-5+j}x_{-3+j}b(1-a^{s-i+1})\neq (1-a)a^{s-i+1},
\end{align}
for $i=0, \;j\in \{0,1\}$; and $i=1,\;j\in\{2,3\}$.
\subsubsection{Case: $a = -1$}
We have\\
\begin{equation*}\label{sole0}
x_{4n - 5} =
\begin{cases}
 c^{1 - n}g^{n}\left(\frac{ -1 + bce }{ -1 + beg }\right)^{\lfloor \frac{n - 1}{2} \rfloor}, & \text{if}\,\, n\,\, \text{is odd};\\\\
c^{1 - n}g^{n}\left(\frac{ -1 + bce }{ -1 + beg }\right)^{\lfloor \frac{n - 1}{2} \rfloor + 1}, & \text{if}\,\, n\,\, \text{is even};
\end{cases}
\end{equation*}
\vspace{1mm}

\begin{equation*}\label{sole1}
 x_{4n - 4} =
 \begin{cases}
d^{1 - n}h^{n}\left(\frac{-1 + bdf}{-1 + bfh}\right)^{\lfloor \frac{n - 1}{2} \rfloor},   & \text{if}\,\, n\,\, \text{is even};\\\\
d^{1 - n}h^{n}\left(\frac{-1 + bdf}{-1 + bfh}\right)^{\lfloor \frac{n - 1}{2} \rfloor + 1},   & \text{if}\,\, n\,\, \text{is odd};
 \end{cases}
\end{equation*}
\vspace{1mm}

\begin{equation*}\label{sole2}
x_{4n - 3} =
\begin{cases}
\frac{c^{n}g^{-n}e}{ -1 + bce }\left(\frac{ -1 + beg}{-1 + bce}\right)^{ \lfloor \frac{n - 1}{2} \rfloor + 1},   & \text{if}\,\, n\,\, \text{is odd};\\\\
c^{n}g^{-n}e\left(\frac{ -1 + beg}{-1 + bce}\right)^{ \lfloor \frac{n - 1}{2} \rfloor + 1},   & \text{if}\,\, n\,\, \text{is even};
\end{cases}
\end{equation*}
\noindent
and
\vspace{0.4mm}
\begin{equation*}\label{sole3}
x_{4n - 2} =
\begin{cases}
\frac{d^{n}h^{-n}f }{ -1 + bdf }\left(\frac{ -1 + bfh}{-1 + bdf}\right)^{ \lfloor \frac{n - 1}{2} \rfloor + 1},   & \text{if}\,\, n\,\, \text{is odd};\\\\
 d^{n}h^{-n}f\left(\frac{ -1 + bfh}{-1 + bdf}\right)^{ \lfloor \frac{n - 1}{2} \rfloor + 1},   & \text{if}\,\, n\,\, \text{is even};
\end{cases}
\end{equation*}
Here, condition \eqref{welldef} simplifies to $bx_{-5+j}x_{-3+j} \neq 1$ for $j = 0,1,2,3$.
\subsection{Case: $a = 1$}
The solution is given by
\begin{equation*}\label{sole0}
x_{4n - 5} =\frac{ g^{n}}{ c^{n-1}}\prod\limits_{s = 0}^{n - 1}\frac{ 1 + bces }{ 1 + begs},\;
x_{4n - 4} =  \frac{h^{n}}{d^{n-1}}\prod\limits_{s = 0}^{n - 1}\frac{1 + bdfs}{1 + bfhs},
\end{equation*}
\begin{equation*}\label{sole2}
x_{4n - 3} =\frac{ c^{n}e}{g^{n}}\prod\limits_{s = 0}^{n - 1}\frac{ 1 + begs }{1 + bce(s+1)} ,\;
x_{4n - 2} = \frac{d^{n}f}{h^{n}}\prod\limits_{s = 0}^{n - 1}\frac{1 + bfhs }{1  + bdf(s + 1)}.
\end{equation*}
Here, condition \eqref{welldef} simplifies to $-x_{-5+j}x_{-3+j} b(s-i+1)\neq 1$ for 
 $i=0, \;j\in \{0,1\}$; and $i=1,\;j\in\{2,3\}$.
\section{The case $a_n$ and $b_n$ are 2-periodic}
In this case, we have $\{a_n\}_{n = 0}^{\infty} = a_0, a_1, a_0, a_1, \ldots$, and similarly $\{b_n\}_{n = 0}^{\infty} = b_0, b_1, b_0, b_1, \ldots$ where $a_0 \neq a_1$, and $b_0 \neq b_1$. Then the solution is given by
\begin{equation*}
x_{4n - 5} = \frac{g^{n}}{c^{ n-1}}\prod\limits_{s = 0}^{n - 1}\frac{a_{0}^{s} + b_0ce\sum\limits_{l = 0}^{s - 1}a_{0}^{l}}{a_{0}^{s} + b_0eg\sum\limits_{l = 0}^{s - 1}a_{0}^{l}},\;
x_{4n - 4} = \frac{h^{n}}{d^{n-1}}\prod\limits_{s = 0}^{n - 1}\frac{a_{1}^{s} + b_1df\sum\limits_{l = 0}^{s - 1}a_{1}^{l}}{ a_{1}^{s} + b_1fh\sum\limits_{l = 0}^{s - 1}a_{1}^{l}},
\end{equation*}
\begin{equation*}
x_{4n - 3} = \frac{c^{n}e}{g^{n}}\prod\limits_{s = 0}^{n - 1}\frac{ a_{0}^{s} + b_0eg\sum\limits_{l = 0}^{s - 1}a_{0}^{l} }{a_{0}^{s + 1} + b_0ce\sum\limits_{l = 0}^{s}a_{0}^{l}},\;
x_{4n - 2} = \frac{d^{n}f}{h^{n}}\prod\limits_{s = 0}^{n - 1}\frac{ a_{1}^{s} + b_1fh\sum\limits_{l = 0}^{s - 1}a_{1}^{l} }{a_{1}^{s + 1} + b_1df\sum\limits_{l = 0}^{s}a_{1}^{l} }.
\end{equation*}
Here, condition \eqref{welldef} simplifies to $-x_{-5+j}x_{-3+j}b_j\sum\limits_{l=0}^{s-i} a_j^l\neq a_j^{s-i+1}$ for
 $i\in\{0,1\}$ and $ \;j\in \{0,1\}$.
\section{The case $a_n, b_n$ are 4-periodic}
We assume that $\{a_n \} = a_0, a_1, a_2, a_3, a_0, a_1, a_2, a_3, \ldots$ and $\{b_n \} =$ $ b_0, b_1, b_2, b_3,$ \\$ b_0, b_1, b_2, b_3, \ldots$. The solution is given by
\begin{equation*}\label{eqn0}
x_{4n - 5} =\frac{ g^{n}}{c^{ n-1}} \prod\limits_{s = 0}^{n - 1}\frac{ a_{0}^{s} + b_0ce\sum\limits_{l = 0}^{s - 1}a_{0}^{l} }{ a_{2}^{s}  + b_2eg\sum\limits_{l = 0}^{s - 1}a_{2}^{l} },\;
x_{4n - 4} = \frac{h^{n}}{d^{ n-1}}\prod\limits_{s = 0}^{n - 1}\frac{ a_{1}^{s} + b_1df\sum\limits_{l = 0}^{s - 1} a_{1}^{l} }{a_{3}^{s}  + b_3fh\sum\limits_{l = 0}^{s - 1}a_{3}^{l} },\;
\end{equation*}
\begin{equation*}\label{eqn2}
x_{4n - 3} = \frac{c^{n}e}{g^{n}}\prod\limits_{s = 0}^{n - 1}\frac{a_{2}^{s} + b_2eg\sum\limits_{l = 0}^{s - 1}a_{2}^{l} }{ a_{0}^{s+1} + b_0ce\sum\limits_{l = 0}^{s} a_{0}^{l}},\;
x_{4n - 2} = \frac{d^{n}f}{h^{n}}\prod\limits_{s = 0}^{n - 1}\frac{a_{3}^{s} + b_3fh\sum\limits_{l = 0}^{s - 1} a_{3}^{l} }{ a_{1}^{s + 1} + b_1df\sum\limits_{l = 0}^{s}a_{1}^{l}}.
\end{equation*}
Here, condition \eqref{welldef} simplifies to $-x_{-5+j}x_{-3+j}b_j\sum\limits_{l=0}^{s-i} a_j^l\neq a_j^{s-i+1}$ for
$i=0, \;j\in \{0,1\}$; and $i=1,\;j\in\{2,3\}$.
\section{Conclusion}
In this paper, non-trivial symmetries for  difference equations of the form \eqref{1.0} were found. Consequently, the results were used to find formulas for the solutions of the equations \eqref{xn}. Specific cases of the solutions were also discussed.

\end{document}